\documentclass[a4paper]{article}
\usepackage[comma,numbers,square,sort&compress]{natbib}
\renewcommand{\top}{{\mathrm{T}}}
\usepackage{graphicx,amsmath,amssymb}
\usepackage{amsthm,bm}
\theoremstyle{definition} %% plain, definition, remark
\usepackage{subfigure,pgf,tikz,pst-node}
\usetikzlibrary{arrows,shapes,automata,positioning,fit}
\usepackage{bm} %%向量的斜体加粗
\newcommand{\lm}{\lambda}
\newcommand{\e}{\varepsilon}
\newcommand{\R}{\mathbb{R}}
\newcommand{\w}{\omega}
\newtheorem{thm}{Theorem}

\newtheorem{lem}{Lemma}
\newtheorem{ass}{Assumption}

\newtheorem{rem}{Remark}

\newcommand{\pb}{\noindent\textbf{Proof. } }
\newcommand{\pe}{\hfill\rule{4pt}{8pt}}

\def\rm{\mathrm}
\begin{document}

\title{Distributed Optimal Steady-State Regulation for High-Order Multi-Agent Systems with External Disturbances}

\author{Yutao Tang \footnote{This work was supported by National Natural Science Foundation of China under Grants 61503033. Yutao Tang is with the School of Automation, Beijing University of Posts and Telecommunications, Beijing 100876. P.\,R. China. Email: yttang@bupt.edu.cn}}

\date{}

\maketitle

{\noindent\bf Abstract}: In this paper, a distributed optimal steady-state regulation problem is formulated and investigated for heterogeneous linear multi-agent systems subject to external disturbances. We aim to steer this high-order multi-agent network to a prescribed steady-state determined as the optimal solution of a resource allocation problem in a distributed way. To solve this problem, we employ an embedded control design and convert the formulated problem to two simpler subproblems. Then, both state-feedback and output feedback controls are presented under mild assumptions to solve this problem with disturbance rejection. Moreover, we extend these results to the case with only real-time gradient information by high-gain control techniques.  Finally, numerical simulations verify their effectiveness.

{\noindent \bf Keywords}: 	optimal steady-state regulation, high-order dynamics, embedded control, disturbance rejection

\section{Introduction}

Multi-agent coordination has been a hot topic over the last decade due to its wide applications in engineering systems. Particularly, steering agents to a prescribed pattern is widely studied in different circumstances. In this type of problems, a (virtual) leader is often set up to describe this pattern or generate references for other agents to follow, while the leader is usually given as a known autonomous dynamic system with possible unknown states. Many  results embodying this idea have been obtained in both theoretical and engineering sides (see \cite{ren2007distributed,su2012cooperative,tang2015distributed,zhang2015distributed,liang2017cooperative} and references therein).

In this paper, we follow this technical line and consider an optimal steady-state regulation problem for a network of high-order agents.  What differentiates our results from existing leader-following literature is that the prescribed pattern of agents here is not given by an autonomous leader but determined as the optimal solution of some resource allocation problem.  We further assume that each agent only knows a part of the cost function and resource data. Consequently, the steady-state for this multi-agent system can not be determined off-line without cooperation among these agents due to the distributedness of this optimization problem. In this way, conventional asymptotic regulation techniques fail to solve this problem for a high-order multi-agent system.

When all agents are single integrators, this problem reduces to a special kind of distributed optimization problems, which has been widely investigated from a viewpoint of mathematical programming. In fact, many effective algorithms were proposed \cite{wang2015cooperative,barreiro2017distributed, yangtao2017distributed,yang2017distributed,wen2018adaptive,xu2017distributed} to solve this problem under various conditions. Different from these results, the decision variables in this paper are determined as outputs of high-order dynamic systems, which can be deemed as some kind of dynamic constraints. Thus, the solvability of this optimization problem over a high-order multi-agent system by a continuous-time algorithm can be much more challenging than the conventional one.  More interestingly, we also consider external disturbances for agents,  which are possibly unbounded and disrupt the optimal performance. To guarantee the applicability of designed controllers, it is necessary to consider the disturbance rejection problems in multi-agent control design. 

In fact, very few results have been obtained in this topic. Some interesting attempts include \cite{trip2016internal} and our preceding work \cite{tang2018iet}, where the considered multi-agent systems are subject to an optimal steady-state regulation constraint. While this constraint is defined over the input side for power networks in \cite{trip2016internal}, we formulate and solve the distributed optimal steady-state regulation problem in \cite{tang2018iet} for a special class of passive systems.  Since the relative degree of passive systems is no higher than one, this motivates us to consider more general agents' dynamics. Although there are similar results to achieve an optimal consensus \cite{zhang2017distributed,wang2016cyber,tang2018cyb}, this is the first time to our knowledge that such a problem is formulated for general linear multi-agent systems subject to external disturbances with a resource-allocation constraint.

In view of the aforementioned observations, the main contributions of this paper are at least two-fold. 
\begin{itemize}
	\item A distributed optimal steady-state regulation problem is formulated and solved for a network of high-order agents. Upon the available information, both state-feedback and output-feedback controls are proposed to solve this problem.  It can be taken as a distributed version of the conventional asymptotic optimal steady-state optimization and  regulation problem for a single plant \cite{skogestad2000self,jokic2009constrained}. 
	\item Compared with existing distributed resource allocation publications for single integrators, this paper considers general linear multi-agent systems by an embedded control approach and thus includes the single-integrator results as special cases \cite{yi2016initialization, cherukuri2015distributed}. 
\end{itemize}
Furthermore, external disturbances are taken into considerations and rejected by a novel observer-based technique under mild assumptions, while only bounded disturbances can be handled in \cite{wang2016cyber,tang2018iet}.

The remainder of this paper is organized as follows.  In Section \ref{sec:formulation}, we present the formulation part. Then, we convert our problem to two subproblems by embedded controls in Section \ref{sec:conversion}. Main results are presented in Section \ref{sec:main} along with two gradient-based controls. Following that, we extend the above results to the case with only real-time gradients in Section \ref{sec:extension} and present several numerical examples to verify the effectiveness of our designs in Section \ref{sec:simu}. Finally, conclusions are given in Section \ref{sec:con}.

\textsl{Notations:} Let $\R^n$ be the $n$-dimensional Euclidean space. For a vector $x$, $||x||$ denotes its Euclidean norm. ${\bf 1}_N$ (and ${\bf 0}_N$) denotes an $N$-dimensional all-one (and all-zero) column vector. $\mbox{col}(a_1,\,{\dots},\,a_n) = [a_1^\top,\,{\dots},\,a_n^\top]^\top$ for column vectors $a_i\; (i=1,\,{\dots},\,n)$. For a given matrix $A\in \R^{m\times m}$, $\sigma(A)$ denotes its spectrum. For a given smooth function $f\colon \R^n\to \R$, we use $\nabla f_i$ and $\nabla^2 f_i$ to represent its gradient and Hessian matrix. Let  $r_N=\frac{1}{\sqrt{N}}{\bf 1}_N$ and $R_N\in \mathbb{R}^{N \times (N-1)}$ satisfying $R_N^\top r_N={\bf 0}_N$, $R_N^\top R_N=I_{N-1}$ and $R_N R_N^\top=I_{N}-r_N r_N^\top$.

\section{Problem formulation}\label{sec:formulation}
Considering a collection of multi-agent systems as follows:
\begin{align}\label{sys:agent}
\begin{split}
\dot{x}_i&=A_ix_i+B_iu_i+E_i \w_i\\
y_i&=C_ix_i, \quad i=1,\,\dots,\, N
\end{split}
\end{align}
where $x_i\in \R^{n_{i}}$, $u_i\in \R$, and $y_i\in \R$ are its state, input, and output variables of agent $i$. The signal $\w_i\in \R^{q_i}$ is an external disturbance of agent $i$ modeled by 
\begin{align}\label{sys:disturbance}
\dot{\w}_i=S_i\w_i.
\end{align}

Associated with these agents, we endow agent $i$ with a local cost function $f_i \colon \mathbb{R} \to \mathbb{R}$. For the whole multi-agent system, we assume its prescribed steady-state is determined as the optimal solution of the following constrained optimization problem:
\begin{align}\label{opt:main}
\begin{split}
\mbox{minimize}&\quad  \quad f(y)=\sum\nolimits_{i=1}^Nf_i(y_i)\\
\mbox{subject to}&\quad    \sum\nolimits_{i=1}^N y_i=\sum\nolimits_{i=1}^N d_i
\end{split}
\end{align}
where $y\triangleq \mbox{col}(y_1,\,\dots, \,y_N)$, the function $f_i(\cdot)$ and constant $d_i$ are private to agent $i$ and can not be shared with other agents. 

This optimization problem is often called resource allocation and many practical applications can be formulated as the above, e.g. economic dispatch in power systems \cite{yi2016initialization}, flow control in networks \cite{bertsekas1998network}. Coupled with the physical agents, we aim to design proper controllers such that the outputs of these agents asymptotically solve the optimization problem \eqref{opt:main}. In other words, the controllers should regulate the agents' outputs such that the equality constraint is satisfied and optimal performance is achieved, both in an asymptotic manner.

Without loss of generality, we assume the triple $(C_i,\,A_i,\,B_i)$ is minimal and the pair $(E_i,\,S_i)$ observable. If not, we can always define a new minimal triple $(C_i,\,A_i,\,B_i)$ and observable pair $(E_i,\,S_i)$ by removing the undesired modes of systems \eqref{sys:agent} and \eqref{sys:disturbance}, which will in no way affect the designed goal\cite{chen1995linear}. 

Due to the privacy of local cost function $f_i(\cdot)$ and resource data $d_i$, each agent can only have a part of the global optimization problem and no agent can gather enough information to determine the optimal solution by itself. Hence, the steady-state can not be obtained off-line without on-line cooperation and information sharing among these agents. 

We use an undirected graph $\mathcal {G}=(\mathcal {N}, \mathcal {E}, \mathcal{A})$ to describe the information sharing topology with the node set
$\mathcal{N}=\{1,{\dots},N\}$ and the edge set $\mathcal {E}$ \cite{mesbahi2010graph}. An edge $(i,\,j)\in \mathcal{E}$ between nodes $i$ and $j$ means that agent $i$ and agent $j$ can share information with each other. The weighted adjacency matrix $\mathcal{A}=[a_{ij}]\in \mathbb{R}^{N\times N}$ is defined by $a_{ii}=0$ and $a_{ij}=a_{ji}\geq0$ ($a_{ij}>0$ if and only if there is an edge between node $i$ and node $j$). The Laplacian $L=[l_{ij}]\in \mathbb{R}^{N\times N}$ of graph $\mathcal{G}$ is defined as $l_{ii}=\sum_{j\neq i}a_{ij}$ and $l_{ij}=-a_{ij} (j\neq i)$, which is thus symmetric. 

Regarding multi-agent system \eqref{sys:agent}, external disturbance \eqref{sys:disturbance} and an information sharing graph $\mathcal{G}$, the distributed optimal steady-state regulation problem is to find distributed controllers for each agent by using its local data and exchanged information with its neighbors, such that its output $y_i(t)$ converges to a steady-state $y_i^*$ as $t$ goes to infinity for $i=1,\,\dots,\,N$, where $\mbox{col}(y_1^*,\,\dots,\,y_N^*)$ is an optimal solution of problem \eqref{opt:main}.

\begin{rem}\label{rem:formulation-2}
	This problem can be taken as a distributed resource allocation problem for heterogeneous high-order agents subject to external disturbances, and thus includes the formulation for single integrators as a special case \cite{yi2016initialization, cherukuri2015distributed}. Particularly, when $f_i(y_i)=\frac{1}{2}y_i^2$ and $d_i=y_i(0)$, the optimal solution of \eqref{opt:main} is $\frac{\sum_{i=1}^N y_i(0)}{N}{\bf 1}_N$. Thus, this formulation provides another way to achieve an average consensus \cite{rezaee2015average} for general liner agents subject to external disturbances.
\end{rem}

To solve the formulated problem, we need several technical assumptions, which are very standard in coordination of multi-agent systems \cite{cherukuri2015distributed,yi2016initialization, tang2015distributed,jakovetic2015linear}.
\begin{ass}\label{ass:graph}
	The   graph $\mathcal{G}$ is  connected.
\end{ass}

\begin{ass}\label{ass:function}
	%	For $i=1,\,\dots,\, N$, the function $f_i\colon\R\to\R$ is twice continuously differentiable with bounded Hessian, i.e., there exist $0 < \underline{h}_i  \leq \bar h_i <\infty$  such that, for all $i$:
	$$\underline{h}_i 	 \leq \nabla^2 f_i(s) \leq \bar h_i ,\quad \forall s \in  \R.$$
\end{ass}

It is well-known that under Assumption \ref{ass:graph}, the associated Laplacian $L$ of this graph is positive semidefinite with rank $N-1$, and its null space is spanned by $\bm{1}_N$. Assumption \ref{ass:function} is made to guarantee the wellposedness of optimization problem \eqref{opt:main}.  In fact, it implies strong smoothness and strong convexity of $f_i(\cdot)$. Then, the problem \eqref{opt:main} is solvable and has a unique solution $ y^*=\mbox{col}(y_1^*,\,\dots, \, y_N^*)$ by the refined Slater's condition \cite{boyd2004convex}. As usual, we assume $y^*$ is finite \cite{yi2016initialization,jakovetic2015linear}.  

\begin{ass}\label{ass:regulator}
	For each $i=1,\,\dots,\,N$, there exist constant matrices $X_{i1}$,\,$X_{i2}$,\,$U_{i1}$ and $U_{i2}$ satisfying
	\begin{align*}
	\begin{split}
	X_{i1}S_i&=A_iX_{i1}+B_iU_{i1}+E_i\\
	0&=C_iX_{i1}
	\end{split}
	\end{align*}
	and
	\begin{align*}
	\begin{split}
	0&=A_iX_{i2}+B_iU_{i2}\\
	1&=C_iX_{i2}.
	\end{split}
	\end{align*}
\end{ass}
\begin{rem}
	This assumption is known as the solvability of regulator equations and plays a key role in achieving the asymptotic steady-state regulation goal \cite{khalil2002nonlinear,huang2004nonlinear}.  A sufficient condition to guarantee this assumption is that
	\begin{align*}
	\mathrm{rank}\begin{bmatrix}
	A_i-\lambda I_{n_i}&    B_i\\
	C_i          &   {\bm 0}
	\end{bmatrix}=n_i+1,\, \forall \lambda\in \sigma(S_i)\cup 0,\,  i=1,\,\dots,\,N.
	\end{align*}
	When agent $i$ satisfies the matching condition, i.e., $B_iM_i=E_i$ for some matrix $M_i$, the first equation naturally holds by letting $X_{i1}=0$ and $U_{i1}=-M_i$. Moreover, for minimum-phase agents,  the second equation is also fulfilled, including integrators as special cases.
\end{rem}

\section{Embedded Control and Problem Conversion}\label{sec:conversion}
To avoid the difficulties brought by the high-order structure of agents, we employ the embedded scheme proposed in \cite{tang2018cyb} to solve the formulated distributed optimal steady-state regulation problem.

In this embedded control scheme, we first introduce an optimal signal generator by considering the same optimization problem for ``virtual" single integrators, in order to asymptotically reproduce the optimal solution $y^*$ by a vector $z\triangleq \mbox{col}(z_1,\,\dots,\,z_N)$. Then, by taking $z_i$ as an output reference signal for agent $i$, we embed this generator in the feedback loop via a tracking controller for system \eqref{sys:agent}. In this way, we can divide the distributed optimal steady-state regulation problem into two simpler subproblems, i.e., optimal signal generator construction and  asymptotic tracker design. The former is to solve the distributed resource allocation problem for single integrators and has been studied by many authors \cite{cherukuri2015distributed,yi2016initialization}, while the latter is an asymptotic tracking problem for general linear agents of the form \eqref{sys:agent} subject to external disturbance generated by \eqref{sys:disturbance}.

To give a self-contained design, we consider the following algorithm proposed in \cite{tang2016distributed}: 
\begin{align}\label{sys:opt-solver}
\begin{split}
\dot{z}_i&=-\nabla f_i(z_i)+\lambda_i\\
\dot{\lambda}_i&=-\lambda^0_i-v^0_i+d_i-z_i\\
\dot{v}_i&=\lambda^0_i
\end{split}
\end{align}
where $\lambda^0_i\triangleq \sum\nolimits_{j=1}^Na_{ij}(\lambda_i-\lambda_j)$, $v^0_i\triangleq\sum\nolimits_{j=1}^N a_{ij}(v_i-v_j)$. 

This algorithm is inherently a primal-dual dynamics to solve the problem \eqref{opt:main} and has also been partially investigated in \cite{yi2016initialization}. The following lemma shows the feasibility of this dynamics as a distributed optimal signal generator for our problem.
\begin{lem}\label{lem:solver-exp}
	Suppose Assumptions \ref{ass:graph}--\ref{ass:function} hold. Then, the algorithm \eqref{sys:opt-solver} will exponentially generate the optimal solution of the distributed optimization problem \eqref{opt:main}, i.e., $z_i$ converges to $y^*_i$ exponentially fast as $t\to+\infty$ for $i=1,\,\ldots,\,N$. %, where $(y_1^*,\,\dots,\, y_N^*)$ is the optimal solution of \eqref{opt:main}.
\end{lem}

\pb First, we show that at the equilibrium point $(z_i^*,\, \lm_i^*,\, v_i^*)$ of system \eqref{sys:opt-solver}, it solves the optimization problem \eqref{opt:main}. For this purpose, we let the righthand side of \eqref{sys:opt-solver} be zero and have $-\nabla f_i(z_i)+\lambda_i=0$,\, $-\lambda^0_i-v^0_i+d_i-z_i=0$,\, $\lambda^0_i=0$. By summing the second equation up from $1$ to $N$,  it follows $-\sum \nolimits_{j=1}^Nv^0_i+\sum \nolimits_{j=1}^N d_i-\sum\nolimits_{j=1}^N z_i=0$. Under Assumption \ref{ass:graph}, we have ${\bm 1}_N^\top L=0$, which implies that $\sum \nolimits_{j=1}^Nv^0_i={\bm 1}_N^\top Lv=0$ and thus $\sum \nolimits_{j=1}^N d_i-\sum\nolimits_{j=1}^N z_i=0$.  Moreover, from $\lambda^0_i=0$ ($i=1,\,\dots,\,N$), we have $L\lambda=\mbox{col}(\lambda^0_1,\,\dots,\,\lambda^0_N)={\bm 0}_N$. This implies that there exists a constant $\theta^*$ satisfying $\lambda_1=\dots=\lambda_N=\theta^*$ since the kernel of graph Laplacian $L$ is spanned by ${\bm 1}_N$. The above arguments imply that at the equilibrium point $(z_i^*,\, \lm_i^*,\, v_i^*)$, it must hold that $\sum \nolimits_{j=1}^N d_i-\sum\nolimits_{j=1}^N z^*_i=0$ and $\nabla f_1(z_1^*)=\dots=\nabla f_N(z_N^*)$.  Then, we recall the refined Slater's condition \cite{boyd2004convex} and conclude that $\mbox{col}(z_1^*,\,\dots,\,z_N^*)$ indeed solves the convex optimization problem \eqref{opt:main}. From the uniqueness of this optimal solution to \eqref{opt:main}, one can further have $y_i^*=z_i^*$ for $i=1,\,\dots,\,N$. 

Next, we consider the exponential stability of \eqref{sys:opt-solver} at the equilibrium $(z_i^*,\, \lm_i^*,\, v_i^*)$. Performing a coordinated transformation $\bar z_i=z_i-z_i^*$, $\bar \lm_i=\lm_i-\lm_i^*$, $\bar v_i=v_i-v_i^*$ gives that
\begin{align*}
\begin{split}
\dot{\bar z}_i&=-{\bm h}_i(\bar z_i,\,z_i^*)+\bar \lambda_i\\
\dot{\bar \lambda}_i&=-\bar \lambda^0_i-\bar v^0_i-\bar z_i\\
\dot{\bar v}_i&=\bar \lambda^0_i
\end{split}
\end{align*}
where ${\bm h}_i(\bar z_i,\,z_i^*)\triangleq \nabla f_i(z_i)-\nabla f_i(z_i^*) $, $\bar \lambda^0_i\triangleq \sum\nolimits_{j=1}^Na_{ij}(\bar \lambda_i-\bar \lambda_j)$, $\bar v^0_i\triangleq\sum\nolimits_{j=1}^N a_{ij}(\bar v_i-\bar v_j)$. Denote $\bar z=\mbox{col}(\bar z_1,\,\dots,\,\bar z_N)$, $\bar \lm=\mbox{col}(\bar \lm_1,\,\dots,\,\bar \lm_N)$ and $\bar v=\mbox{col}(\bar v_1,\,\dots,\,\bar v_N)$.  It can be written into a compact form: 
\begin{align*}
\begin{split}
\dot{\bar z}&=-{\bm h}(\bar z,\,z^*)+\bar \lambda\\
\dot{\bar \lambda}&=-L\bar \lambda-L\bar v-\bar z\\
\dot{\bar v}&=L\bar \lambda
\end{split}
\end{align*}
where the vector-valued function ${\bm h}(\bar z,\,z^*)$ is determined by ${\bm h}_i(\bar z_i,\,z_i^*)$. By Assumption \ref{ass:function}, ${\bm h}(\bar z,\,z^*)$ is $\bar h$-Lipschitz in $\bar z$ and satisfies $\bar z^\top {\bm h}(\bar z,\,z^*)\geq \underline{h} \bar z^\top \bar z$ for $\overline{h}=\max_i\{\overline{h}_i \}$ and $\underline{h}=\min_i\{\underline{h}_i \}$. 

At present, we only have to investigate the stability of the above system at the origin. For this purpose, let $\hat v_1=r^\top \bar v$, $\hat v_2=R^\top \bar v$. Since $r^\top \dot{v}=0$, it follows that $\dot{\hat v}_1=r^\top \dot{\bar v}=0$. This means $r^\top v(t)$ is an invariant quantity and thus has no effect on the stability of the above system. Furthermore, the rest states of this system can be put into the following form: 
\begin{align*}
\begin{split}
\dot{\bar z}&=-{\bm h}(\bar z,\,z^*)+\bar \lambda\\
\dot{\bar \lambda}&=-L\bar \lambda-LR\hat v_2-\bar z\\
\dot{\hat v}_2&=R^\top L^\top \bar \lambda
\end{split}
\end{align*}
where we use $Lr=0$ and $L=L^\top$ by Assumption \ref{ass:graph}.   

Let $T\triangleq \begin{bmatrix}-L &-LR \\ R^\top L^\top &{\bm 0}\end{bmatrix}$ and $S\triangleq \begin{bmatrix}
-I_{N} \\ {\bm 0}
\end{bmatrix}^\top $. It can be easily verified that $T+T^\top \leq {\bm 0}$ and the pair $(S,\,T)$ is observable. Then, we apply Lemma 2 in \cite{tang2018cyb} to the above system and conclude its exponential stability at the origin, which implies the conclusions. The proof is thus complete.
\pe

With the designed optimal signal generator, our formulated distributed optimal steady-state regulation problem is converted to an asymptotic tracking and disturbance rejection problem for linear multi-agent systems of the form \eqref{sys:agent} subject to external disturbances generated by \eqref{sys:disturbance}. %Next, we will propose two kinds of tracking controllers for the considered agents and finally solve our problem. 

\section{Solvability of Optimal Steady-State Regulation}\label{sec:main}

In this section, we seek  proper tracking controllers for agent \eqref{sys:agent} with output reference $z_i$, and finally solve the associated optimal steady-state regulation problem. 

Note that $(A_i,\, B_i)$ is controllable, there exists $K_{i1}$ such that the matrix $A_i+B_iK_{i1}$ is Hurwitz. Denote $K_{i2}\triangleq U_{i1}-K_{i1}X_{i1}$, $K_{i3}\triangleq U_{i2}-K_{i1}X_{i2}$. From the conventional output regulation theory (e.g., Theorem 1.7 in \cite{huang2004nonlinear}), a full-information regulator $u^o_i=K_{i1}x_i+K_{i2}\w_i+K_{i3}y_i^*$ can steer agent $i$ to its associated steady-state and also achieve disturbance rejection.  As the optimal signal generator \eqref{sys:opt-solver} can exponentially generate $y_i^*$, we only have to consider the disturbance rejection part. 

Inspired by the techniques used in \cite{huang2004nonlinear,tang2016dobc},  we construct a reduced-order observer to estimate the disturbance $\w_i$.
\begin{align}\label{sys:reduced-observer}
\begin{split}
\dot{\bar \eta}_i&=(S_i-\bar L_i E_i)\bar \eta_i+(S_i\bar L_i-\bar L_iE_i\bar L_i-\bar L_iA_i)x_i-\bar L_iB_iu_i\\
\eta_i&=\bar \eta_i+\bar L_i x_i
\end{split}
\end{align}
where $\bar L_i$ is chosen such that the matrix $S_i-\bar L_iE_i$ is Hurwitz. Its existence is guaranteed by the observability of $(E_i,\,S_i)$.  

Here is a lemma to show the effectiveness of this estimator.
\begin{lem}
	Consider the dynamic system composed of \eqref{sys:agent} and \eqref{sys:reduced-observer}. Then, $\eta_i(t)$ exponentially converges to $\w_i$, i.e., there exist two positive constants $\mu_1,\, \mu_2$ satisfying $||\eta_i(t)-\w_i(t)||\leq \mu_1e^{-\mu_2 t}$ for all $t>0$.
\end{lem}
\pb 
To prove this lemma, we let $\hat \eta_i\triangleq\eta_i-\w_i$. Then, along the trajectory of systems \eqref{sys:disturbance} and \eqref{sys:reduced-observer}, we have 
\begin{align*}
\dot{\hat \eta}_i&=\dot{\overline \eta}_i+ \bar L_i\dot{x}_i-\dot{\w}_i\\
&=(S_i-\bar L_i E_i)\bar \eta_i+(S_i\bar L_i-\bar L_iE_i\bar L_i-\bar L_iA_i)x_i\\
&-\bar L_iB_iu_i+\bar L_i(A_ix_i+B_iu_i+E_i \w_i)-S_i\w_i\\
&=(S_i-\bar L_i E_i)\hat \eta_i.
\end{align*}
By the selection of $\bar L_i$, the matrix $(S_i-\bar L_i E_i)$ is Hurwitz. In other words, this error system is exponentially stable. This implies the existence of $\mu_1,\,\mu_2>0$ satisfying this lemma.  \pe

With this estimator, we substitute it into the full-information regulator and propose a state-feedback controller for agent \eqref{sys:agent} to solve our formulated problem as follows.

\begin{thm}\label{thm:main}
	Suppose Assumptions \ref{ass:graph}--\ref{ass:regulator} hold.  Then, the distributed optimal steady-state regulation problem for multi-agent system \eqref{sys:agent} with optimization problem \eqref{opt:main} and disturbance \eqref{sys:disturbance} is exponentially solved by the following control:
	\begin{align}\label{ctrl:main}
	\begin{split}
	u_i&=K_{i1}x_i+K_{i2}\eta_i+K_{i3}z_i\\
	\dot{\bar \eta}_i&=(S_i-\bar L_i E_i)\bar \eta_i+(S_i\bar L_i-\bar L_iE_i\bar L_i-\bar L_iA_i)x_i-\bar L_iB_iu_i\\
	\dot{z}_i&=-\nabla f_i(z_i)+\lambda_i\\
	\dot{\lambda}_i&=-\lambda^0_i-v^0_i+d_i-z_i\\
	\dot{v}_i&=\lambda^0_i
	\end{split}
	\end{align}
	where $\eta_i=\bar \eta_i+\bar L_i x_i$, $K_{i1},\,K_{i2},\,K_{i3}$ and $\bar L_i$ are chosen gain matrices as above.
\end{thm}

\pb To prove this theorem, we first perform a coordination transformation $\bar x_i=x_i-X_{i1}\w_i- X_{i2} z_i$ and $\hat \eta_i=\eta_i-\w_i$. By Assumption \ref{ass:regulator} and some mathematical manipulations, the whole composite system can be put into the following form:
\begin{align}\label{sys:composite-thm:main}
\begin{split}
\dot{\bar x}_i&=(A_i+B_iK_{i1})\bar x_i+B_iK_{i2}\hat \eta_i-X_{i2}\dot{z}_i\\
\dot{\hat \eta}_i&=(S_i-\bar L_i E_i)\hat \eta_i \\
\dot{z}_i&=-\nabla f_i(z_i)+\lambda_i\\
\dot{\lambda}_i&=-\lambda^0_i-v^0_i+d_i-z_i\\
\dot{v}_i&=\lambda^0_i\\
\bar y_i&=C_i \bar x_i
\end{split}
\end{align}
where $\bar y_i\triangleq y_i-z_i$. 

By the selection of $K_{i1}$, the matrix $\bar A_i \triangleq A_i+B_i K_{i1}$ is Hurwitz. Thus, the Lyapunov equation $\bar A_i^\top P_i+P_i\bar A_i=-2I_{n_i}$ has a unique positive definite solution $P_i$. Similarly, we can determine a unique positive definite matrix $Q_i$ satisfying $(S_i-\bar L_i E_i)^\top Q_i+Q_i (S_i-\bar L_i E_i)=-I_{q_i}$. 

Considering the first two subsystems, we let $V^0_i(\bar x_i,\,\hat \eta_i)=\bar x_i^\top P_i\bar x_i+ c_i \hat \eta_i^\top Q_i\hat \eta_i$ with a constant $c_i>0$ to be specified later. Its time derivative along the trajectory of \eqref{sys:composite-thm:main} is:
\begin{align*}
\dot{V}^0_i&=2\bar x_i^\top P_i[\bar A_i \bar x_i+B_iK_{i2}\hat \eta_i-X_{i2}\dot{z}_i]+2c_i\hat \eta_i^\top Q_i (S_i-\bar L_i E_i)\hat \eta_i\\
&=-2\bar x_i^\top \bar x_i + 2\bar x_i^\top P_i B_iK_{i2}\hat \eta_i-2\bar x_i^\top P_i X_{i2}\dot{z}_i-c_i\hat \eta_i^\top \hat \eta_i.
\end{align*}
By Young's inequality, it follows that 
\begin{align*}
\dot{V}^0_i&\leq -\bar x_i^\top\bar x_i + 2||P_i B_iK_{i2}||^2||\hat \eta_i||^2+2||P_i X_{i2}||^2||\dot{z}_i||^2-c_i\hat \eta_i^\top \hat \eta_i.
\end{align*}
Letting $c_i>2||P_i B_iK_{i2}||^2+1$ implies 
\begin{align*}
\dot{V}^0_i&\leq -\bar x_i^\top \bar x_i -||\hat \eta_i||^2+2||P_i X_{i2}||^2||\dot{z}_i||^2.
\end{align*}
We let $V^0=\sum\nolimits_{i=1}^NV_i^0$ and further have
\begin{align*}
\dot{V}^0&\leq \sum\nolimits_{i=1}^N[-\bar x_i^\top \bar x_i -||\hat \eta_i||^2+2||P_i X_{i2}||^2||\dot{z}_i||^2]\\
&\leq -||\bar x||^2 -||\hat \eta||^2 + 2\sum\nolimits_{i=1}^N ||P_i X_{i2}||^2||\dot{z}_i||^2.
\end{align*}
Choosing $c=1+2\max\limits_i\{||P_i X_{i2}||^2\}$ gives
\begin{align*}
\dot{V}^0\leq -||\bar x||^2 -||\hat \eta||^2 + c ||\dot{z}||^2.
\end{align*}

Next, we consider the optimal signal generator part. From the proof of Lemma \ref{lem:solver-exp}, the following system is exponentially stable at the origin.
\begin{align*}
\begin{split}
\dot{\bar z}&=-{\bm h}(\bar z,\,z^*)+\bar \lambda\\
\dot{\bar \lambda}&=-L\bar \lambda-LR\hat v_2-\bar z\\
\dot{\hat v}_2&=R^\top L^\top \bar \lambda
\end{split}
\end{align*}
Denote $\hat z\triangleq\mbox{col}(\bar z,\, \bar \lambda,\, {\hat v}_2)$ for short.  Since ${\bm h}(\bar z,\,z^*)$ is global Lipschitz for some $\hat l>1$ in $\bar z$ by Assumption \ref{ass:function}, invoking the converse Lyapunov theorem (Theorem 4.15 in \cite{khalil2002nonlinear}) provides us a Lyapunov function $W(\hat z)$ satisfying the following conditions for some positive constants $\hat c_1$,\,$\hat  c_2$,\,$\hat c_3$,\,$\hat c_4$:
\begin{align*}
&\hat c_{1}||\hat z||^2\leq W(\hat z)\leq \hat c_{2}||\hat z||^2\\
&\dot{W}\leq -\hat c_{3}||\hat z||^2,\quad ||\frac{\partial {W}}{\partial{\hat z}}|| \leq \hat c_{4}||\hat z||.
\end{align*}
Since $\dot{z}=\dot{\bar z}$, we choose a Lyapunov function for the composite system as $V(\bar x,\,\hat \eta,\,\hat z)=V^0(\bar x,\, \hat \eta)+\hat c W(\hat z)$.  Its time derivative along the trajectory of the closed-loop system satisfies
\begin{align*}
\dot{V}&\leq -||\bar x||^2 -||\hat \eta||^2 + c ||\dot{z}||^2-\hat c \hat c_{3}||\hat z||^2\\
&\leq -||\bar x||^2 -||\hat \eta||^2 + \frac{c\hat l^2}{\hat c_1} W(\hat z) -\frac{\hat c \hat c_{3}}{\hat c_2}W(\hat z)\\
&\leq -||\bar x||^2 -||\hat \eta||^2 -(\frac{\hat c \hat c_{3}}{\hat c_2}- \frac{c\hat l^2}{\hat c_1}) W(\hat z)
\end{align*}
where we use the $\hat l$-Lipschitzness of $\dot{z}$. 

Letting $\hat c\geq\frac{\hat c_2(\hat c_1 +c\hat l^2)}{\hat c_1\hat c_3}$ gives
\begin{align*}
\dot{V}&\leq -||\bar x||^2 -||\hat \eta||^2 - W(\hat z).
\end{align*}
Recalling the quadratic form of $V^0$ and the properties of $W$, we invoke Theorem 4.10 in \cite{khalil2002nonlinear} and conclude the global exponential convergence of the whole composite system with respect to its equilibrium point. 

From the above arguments, we have that a) $z_i$ exponentially converges to $y_i^*$ as $t\to \infty$ and b) $y_i$ exponentially converges to $z_i$ as $t\to \infty$.  It follows by combing these two facts that
\begin{align*}
|y_i-y_i^*|\leq |y_i-z_i|+|z_i-y_i^*|\to 0 \textrm{ as }  t \to \infty.
\end{align*}
Moreover, the convergence is also exponentially fast.  The proof is thus complete.
\pe

\begin{rem}\label{rem:optimalizing}
	Compared with existing steady-state regulation problem \cite{skogestad2000self,jokic2009constrained}, we consider its distributed extension for a multi-agent system where the steady-state can only be determined and reached in a distributed way, which is of course more challenging. When agents are all single integrators, the above results are consistent with existing distributed resource allocation conclusions \cite{cherukuri2015distributed,yi2016initialization}. Even in this case, we consider unknown (unbounded) disturbances, while most of relevant publications are disturbance-free.
\end{rem}

In some cases, only the output variables of each agent can be obtained because it may be difficult to get or measure all the state variables in some situations. Since the optimal signal generator is independently implemented, we only have to focus on the tracking part to solve the optimal regulation problem with output feedbacks. 

By attaching an observer to the above embedded controller \eqref{ctrl:main}, the following theorem is readily obtained.
\begin{thm}\label{thm:main-output}
	Suppose Assumptions \ref{ass:graph}--\ref{ass:regulator} hold.  Then, the distributed optimal steady-state regulation problem for multi-agent system \eqref{sys:agent} with optimization problem \eqref{opt:main} and disturbance \eqref{sys:disturbance} is exponentially solved by the following control:
	\begin{align}\label{ctrl:main-output}
	\begin{split}
	u_i&=K_{i1}\xi_i+K_{i2}\eta_i+K_{i3}z_i\\
	\dot{\xi}_i&=A_i\xi_i + B_iu_i+ \hat L_i(C_i\xi_i-y_i)\\
	\dot{\bar \eta}_i&=(S_i-\bar L_i E_i)\bar \eta_i+(S_i\bar L_i-\bar L_iE_i\bar L_i-\bar L_iA_i)\xi_i-\bar L_iB_iu_i\\
	\dot{z}_i&=-\nabla f_i(z_i)+\lambda_i\\
	\dot{\lambda}_i&=-\lambda^0_i-v^0_i+d_i-z_i\\
	\dot{v}_i&=\lambda^0_i
	\end{split}
	\end{align}
	where the gain matrices $K_{i1}$,\,$K_{i2}$,\,$K_{i3}$, $\bar L_i$ are defined as in controller \eqref{ctrl:main}, $\hat L_i$ is chosen such that $A_i+\hat L_i C_i$ is Hurwitz.
\end{thm}

The proof is similar to that of Theorem \ref{thm:main} and thus omitted.

\section{Extension with Real-time Gradients}\label{sec:extension}

In many practical cases, we can not have the cost function $f_i(\cdot)$ itself and only its real-time gradient $\nabla f_i(y_i)$ is available. Then, the above control laws will not be implementable. 

Note that there will always be some error between $\nabla f_i(y_i)$ and $\nabla f_i(z_i)$ when $z_i\neq y_i$ by Assumption \ref{ass:function}, and this error will be smaller if the tracking subsystem evolves in a faster time scale. For this purpose, we focus on minimum-phase agents of the form \eqref{sys:agent} with relative degree $r_i>0$, and use a high-gain rule to choose the gain matrices. We only consider the state feedback case here, while its output feedback extension can be derived without difficulties and thus omitted.

Let $\bar X_{i1,1}=0$, $\bar X_{i1,k+1}=\bar X_{i1,k}S_i-C_iA_i^{k-1}E_i$ for $k=1,\,\dots,\,r_i-1$. Denote $\bar X_{i1}=\mbox{col}(\bar X_{i1,1},\,\dots,\,\bar X_{i1,r_i})$, $\bar U_{i1}\triangleq \frac{\bar X_{i1,r_i}S_i-C_iA_i^{r_i-1}E_i}{C_iA_i^{r_i-1}B_i}$, $\bar X_{i2}\triangleq \mbox{col}(1,\,0,\,\cdots,\,0)$,   $U_{i2}\triangleq 0$,  $\bar K_{i1}\triangleq -\frac{1}{C_iA_i^{r_i-1}B_i\e^{r_i}}[c_{i0}~\e c_{i1}~\cdots~\e^{r_i-1}c_{ir_i}]$, $\bar K_{i2}\triangleq \bar U_{i1}-\bar K_{i1}\bar X_{i1}$, $K_{i3}\triangleq \bar U_{i2}-\bar K_{i1}\bar X_{i2}$, and $\hat X_{i}\triangleq \mbox{col}(C_i,\,C_iA_i,\,\dots,\,C_iA_i^{r_i-1})$ with constants $\e,\,c_{i0},\,\dots,\,c_{ir_i-1}$ to be specified later.

The main result with real-time gradients is as follows.
\begin{thm}\label{thm:main-online}
	Suppose Assumptions \ref{ass:graph}--\ref{ass:function} hold and agent \eqref{sys:agent} is minimum-phase. Then, there exists a constant $\varepsilon^*>0$ such that the distributed optimal steady-state regulation problem for multi-agent system \eqref{sys:agent} with optimization problem \eqref{opt:main} and disturbance \eqref{sys:disturbance} is exponentially solved by the following control with any $\varepsilon\in(0, \,\varepsilon^*)$.
	\begin{align} \label{ctrl:main-online}
	u_i=&-\frac{C_i A_i^{r_i}}{C_iA_i^{r_i-1}B_i}x_i+\bar K_{i1}\hat X_{i}x_i+\bar K_{i2}\bar \eta_i+\bar K_{i3}z_i \nonumber \\
	\dot{\bar \eta}_i=&(S_i-\bar L_i E_i)\bar \eta_i+(S_i\bar L_i-\bar L_iE_i\bar L_i-\bar L_iA_i)x_i-\bar L_iB_iu_i\nonumber \\
	\dot{z}_i=&-\nabla f_i(y_i)+\lambda_i\\
	\dot{\lambda}_i=&-\lambda^0_i-v^0_i+d_i-z_i\nonumber \\
	\dot{v}_i=&\lambda^0_i\nonumber
	\end{align}
	where the matrix $\bar L_i$ is defined as above and the constants $c_{i0},\dots,\,c_{ir_i-1}$ are chosen such that the polynomial $\sum_{k=0}^{r_i-1}c_{ik}s^{k}+s^{r_i}$ is Hurwitz, $i=1,\,\dots,\, N$.
\end{thm}
\pb First, we put the original dynamics into a normal form. Note that the plant \eqref{sys:agent} has a relative degree $r_i$. When $r_i\geq 2$,  we can always define a coordinate transformation and put system \eqref{sys:agent} into the following form:
\begin{align*}
\dot{\chi}_{i1}&={\chi_{i2}}+C_iE_i\w_i\\
\vdots\\
\dot{\chi}_{ir_i}&=C_iA_i^{r_i}x_i+C_iA_i^{r_i-1}B_iu_i+C_iA_i^{r_i-1}E_i\w_i\\
\dot{\chi}^z_{i}&=A_i^0\chi^z_{i}+B_i^z\chi_{i}+E_i^z\w_i\\
y_i&=\chi_{i1}
\end{align*}
where $\chi_{ik}=C_iA_i^{k-1}x_i$ and $\chi_i^z$ is the rest state variable. When $r_i=1$, the transformed system is described by a similar form: 
\begin{align*}
\dot{\chi}_{i1}&=C_iA_ix_i+C_iB_iu_i+C_iE_i\w_i\\
\dot{\chi}^z_{i}&=A_i^0\chi^z_{i}+B_i^z\chi_{i}+E_i^z\w_i\\
y_i&=\chi_{i1}.
\end{align*}

To simplify the following proof, we only consider the case when $r_i\geq 2$ without loss of generality. The arguments for $r_i=1$ is almost the same and thus omitted. 

From the minimum-phase assumption, matrix $A_i^0$ is Hurwitz, which implies that the trajectory of $\chi^z_i$ is always well-defined.  By letting $u_i=\bar u_i-\frac{C_i A_i^{r_i}}{C_iA_i^{r_i-1}B_i}x_i$, the $\chi^z_{i}$-subsystem will not effect the output of this plant. Thus, we only have to consider the $\chi_i$-subsystem for agent $i$ with this new input as follows:
\begin{align*}
\dot{\chi}_{i1}&={\chi_{i2}}+C_iE_i\w_i\\
\vdots\\
\dot{\chi}_{ir_i}&=C_iA_i^{r_i-1}B_i\bar u_i+C_iA_i^{r_i-1}E_i\w_i\\
y_i&=\chi_{i1}.
\end{align*}

Next, we solve the distributed optimal steady-state regulation problem for the above agents with disturbance \eqref{sys:disturbance} and optimization problem \eqref{opt:main}. 

Considering the closed-loop system, we substitute the control \eqref{ctrl:main-online} into each agent and define a coordinate transformation $\bar \chi_i=\chi_i-\bar X_{i1}\w_i-\bar X_{i2}z_i$. It follows that
\begin{align*}
\dot{{\bar\chi}}_{i1}&={\bar \chi}_{i2}-\dot{z}_i\\
\vdots\\
\dot{{\bar \chi}}_{ir_i}&=-\frac{1}{\e^{r_i}}[c_{i0}\bar \chi_{i1}+\sum\limits_{k=1}^{r_i-1}\e^{k} c_{ik} \bar \chi_{i(k+1)}]\\
\bar y_i&=\bar \chi_{i1}.
\end{align*}
Letting $\hat \chi_{i1}=\bar \chi_{i1}$ and $\hat \chi_{ik}=\varepsilon^{k-1}\bar \chi_{ik}$ ($k=2,\,\dots,\, r_i$) gives
\begin{align*}
\dot{\hat{\chi_i}}=\frac{1}{\varepsilon} A_{ic} \hat \chi_i-B_{ic} \dot{z}_i
\end{align*}
where $A_{ic}={\left[\begin{array}{c|c}
	0 & \mathrm{I}_{r_i-1}\\ \hline
	c_{i0}&c_{i1}\, \dots\, c_{ir_i-1}
	\end{array}\right]}$ and $B_{ic}=\mbox{col}(1,\,0,\,\dots,\,0)$. Putting all agents in a compact form gives
\begin{align}\label{sys:sp-2}
\varepsilon \dot{\hat{\chi}}=  A_{c} \hat \chi - \varepsilon B_{c} \dot{z}
\end{align}
where $\hat \chi\triangleq \mbox{col}(\hat \chi_{1},\,\dots,\, \hat \chi_N)$,\, $z\triangleq \mbox{col}(z_{1},\,\dots,\, z_N)$,\,  $A_c\triangleq \mbox{block}\,diag \{A_{1c},\, \dots,\, A_{Nc}\}$ and  $B_c\triangleq \mbox{block}\,diag \{B_{1c},\, \dots,\, B_{Nc}\}$. By constructions, $A_{ic}$ and $A_c$ are Hurwitz.

By similar arguments in the proof of Lemma \ref{lem:solver-exp}, the optimal signal generator part can be rewritten as follows:
\begin{align}\label{sys:sp-1}
\begin{split}
\dot{\bar z}&=-{\bm h}(z, y^*)+\Delta(z,\,y)+\bar \lambda-{\bm h}(y,\,z)\\
\dot{\bar \lambda}&=-L\bar \lambda-LR \hat v_2 - \bar z\\
\dot{\hat v}_2&=R^\top L\bar \lambda
\end{split}
\end{align}
where $\bar z=z-y^*$ and $\Delta(z, y)$ is a vector-value function defined as follows:
\begin{align*}
\Delta(z,\,y)\triangleq [\nabla f_1(z_1)-\nabla f_1(y_1),\, \dots,\, \nabla f_N(z_N)-\nabla f_N(y_N)]^\top.
\end{align*}
By Assumption \ref{ass:function},  $\Delta(z, y)$ is global Lipschitz in $\bar y$ and then $\hat \chi$.

Note that when $y=z$,  we have ${\bm h}(y,\,z)=0$ and the above system is exponentially stable at the origin by Lemma \ref{lem:solver-exp}.  Since $\dot{z}=\dot{\bar z}$, it implies that $\dot z$ is linearly upper bounded by ${\rm col}(\bar z,\,\bar \lm,\, \hat \chi)$. Then, the whole system composed of \eqref{sys:sp-2} and \eqref{sys:sp-1} is in a singularly perturbed form. By Corollary 2.3 in Chapter 7 of \cite{kokotovic1999singular}, we can obtain the conclusions.
\pe

\begin{rem}
	The embedded control approach used here was first proposed in \cite{tang2018cyb} to solve an optimal output consensus problem of linear multi-agent systems. Here, we extend it to tackle a much different problem, where the global optimization problem is constrained and the agents are subject to external disturbances. With this approach, the design complexities brought by high-order dynamics and disturbance direction are decoupled from the optimization task. Then, the optimal steady-state regulation problem is solved in a constructive way by solving two simpler subproblems, which verifies the flexibility and effectiveness of this embedded design.
\end{rem}

\section{Simulations}\label{sec:simu}
In this section, we provide some numerical examples to illustrate the effectiveness of previous designs.

{\em Example 1}. Consider an optimal rendezvous problem of wheeled robots with the following dynamics.
\begin{align*}
\begin{cases}
\dot{r}_i^x=v_i\cos(\theta_i)\\
\dot{r}_i^y=v_i\sin(\theta_i)\\
\dot{\theta}_i=\omega_i\\
\dot{v}_i=\frac{1}{m_i}F_i\\
\dot{\omega}_i=\frac{1}{J_i}\tau_i
\end{cases}
\end{align*}
where $(r_i^x,\, r_i^y, \,\theta_i)$ are the inertia center's position and orientation of the $i$th robot, $(v_i,\,\omega_i)$ the linear and angular speed, $(F_i,\,\tau_i)$ the applied force and torque, $(m_i,\, J_i)$ the mass and moment of inertia for $i=1,\,\dots,\,4$. Let $c_i$ represent the distance between the hand position and inertia center of the $i$th robot. Following the arguments in \cite{ren2007distributed}, we employ an output feedback linearization technique about the hand position $\tilde x_i \triangleq (r_i^x+c_i\cos(\theta_i), \,r_i^y+c_i\sin(\theta_i))$ and obtain a simple linear dynamics as follows:
\begin{align*}
\dot{\tilde x}_i=\tilde v_i,\quad \dot{\tilde  v}_i=\tilde u_i,\quad y_i=\tilde x_i.
\end{align*}

To drive all hands of robots to rendezvous at a common point that minimizes the aggregate distance from their starting points to this final location, we formulate this problem as an optimal steady-state regulation problem by taking the cost functions as $f_i(y_i)=\frac{1}{2} ||y_i||^2$ and $d_i=y_i(0)$ ($i=1,\,\dots,\, 4$). We can easily check that the optimal solution of the global cost function is $y^*= \frac{\sum_{i=1}^4y_i(0)}{4}{\bm 1}_4$. 

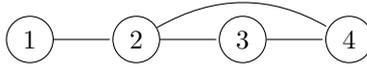
\begin{figure}
	\centering
	\begin{tikzpicture}[shorten >=1pt, node distance=1.4cm, >=stealth',
	every state/.style ={circle, minimum width=0.3cm, minimum height=0.3cm}, auto]
	\node[align=center,state](node1) {1};
	\node[align=center,state](node2)[right of=node1]{2};
	\node[align=center,state](node3)[right of=node2]{3};
	\node[align=center,state](node4)[right of=node3]{4};
	\path[-]   (node1) edge (node2)
	(node2) edge [bend left] (node4)
	(node2) edge  (node3)
	(node3) edge   (node4)
	;
	\end{tikzpicture}
	\caption{The information sharing graph $\mathcal{G}$.}\label{fig:graph}
\end{figure}

The information sharing graph satisfying Assumption \ref{ass:graph} is chosen as Fig.~\ref{fig:graph} with unity edge weights. Choose constants $c_{i0}$ and $c_{i1}$ such that the polynomial $s^2+c_{i1}s+c_{i0}$ is Hurwitz. Then, the optimal rendezvous problem is solved according to Theorem \ref{thm:main} by a control as follows:
\begin{align}\label{ctrl:ex1}
\begin{split}
u_i&=-c_{i0}(y_i-z_i)-c_{i1} \tilde v_i \\
\dot{z}_i&=-y_i+\lambda_i \\
\dot{\lambda}_i&=-\lambda^0_i-v^0_i+d_i-z_i \\
\dot{v}_i&=\lambda^0_i.  
\end{split}
\end{align}

\begin{figure}
	\centering
	\includegraphics[width=0.84\textwidth]{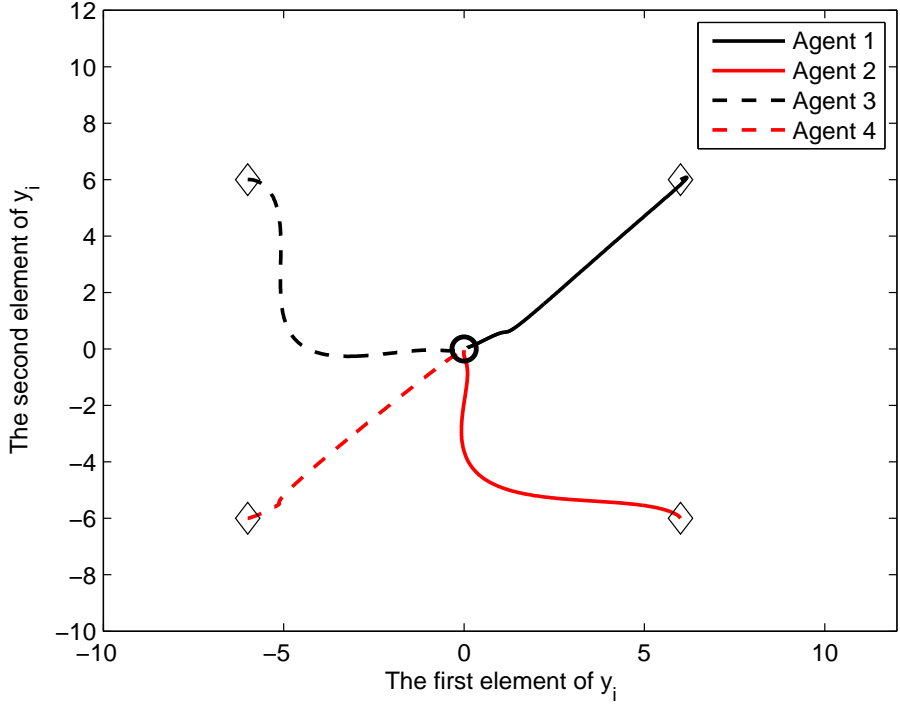}\\
	\caption{Profiles of agents' positions under control \eqref{ctrl:ex1}.}\label{fig:simu-1}
\end{figure}

Take $c_{i0}=4, \, c_{i1}=8,\,\e=1$ with initials (randomly) chosen in $[-10,\,10]$. The simulation result is given in Fig.~\ref{fig:simu-1} and all robots achieve the optimal rendezvous at $\frac{\sum_{i=1}^4y_i(0)}{4}$.

{\em Example 2}. Another example is the distributed inventory control problem with unknown demand rates.  Consider a network of inventories which produce only one commodity \cite{raafat1991survey}. The inventory system at node $i$ can be modeled as
\begin{align}\label{sys:inventory}
\dot{I}_i= P_i-D_i,
\end{align}
where $I_i$, $P_i$ and $D_i$ are the inventory level, production rate and demand rate at node $i$. The storage cost at each warehouse is given as $f_i(I_i)=\alpha_i I_i^2+\beta_iI_i+\gamma_i$, where $\alpha_i>0$.  In distributed inventory control, we aim to maintain the total inventory at certain level $I^{r}$ to satisfy the customer's demands and some safety goals.  To make it more interesting, we consider the case when the demand rate $D_i$ is constant but unknown.

This problem can be formulated as a distributed optimal steady-state regulation problem for inventory systems with external disturbances, while the optimal steady-state is determined by the following problem.  
\begin{align}\label{opt:inventory}
\begin{split}
\mbox{minimize}&\quad  \quad \sum\nolimits_{i=1}^Nf_i(I_i)\\
\mbox{subject to}&\quad   ~~\sum\nolimits_{i=1}^N I_i=I^r.
\end{split}
\end{align}
Clearly, with a preallocation of inventory level $I^r=\sum\nolimits_{i=1}^N I_i^r$, all assumptions are fulfilled to apply the control \eqref{ctrl:main}. 

\begin{figure}
	\centering
	\includegraphics[width=0.84\textwidth]{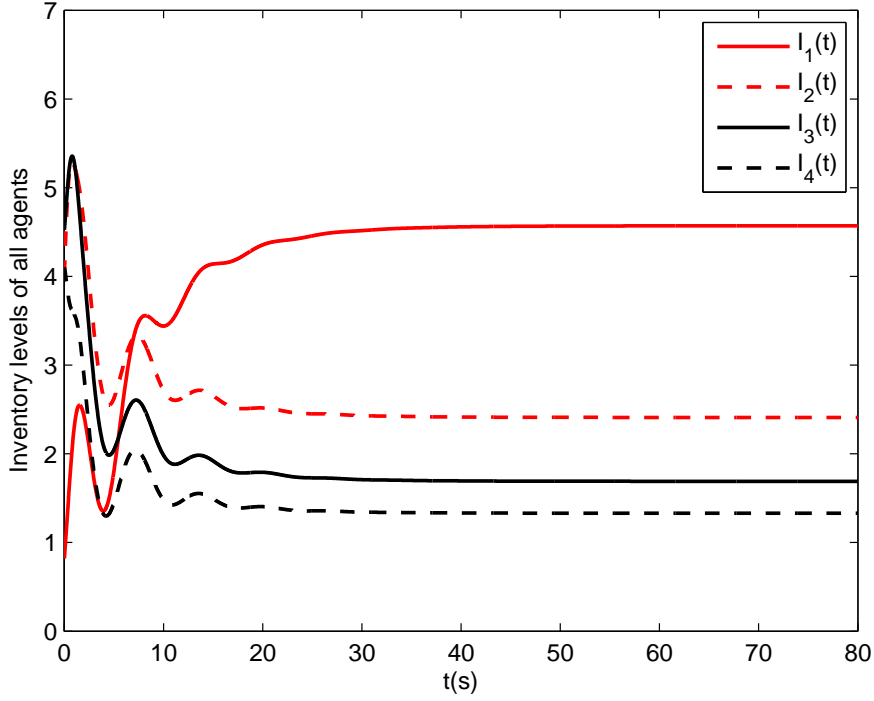}\\
	\caption{Profiles of all inventory levels  under state feedback control \eqref{ctrl:main}.}\label{fig:simu-2}
\end{figure}

For simulations, we take $N=4,\, \alpha_i=0.1 i,\, \beta_i=-0.05 i,\,  \gamma_i=D_i=I_i^r=i$, $i=1,\,\dots,\, 4$. The communication graph is still taken as Fig.~\ref{fig:graph} with initials picked in $[0,\,6]$. By choosing proper gain matrices, we solve this problem and drive the inventory levels of all agents to the optimal solution $I^*=\mbox{col}(4.57,\, 2.41,\,1.69,\, 1.33)$ as depicted in Fig.~\ref{fig:simu-2}.

{\em Example 3}.  Consider a distributed coordination problem for a multi-agent system consisting of four mass-damper-spring systems with unit
mass described by:
\begin{align*}
\ddot{y}_i+g_i\dot{y}_i+f_iy_i=u_i+F_i\w_i,\quad i=1,\,2,\,3,\,4
\end{align*}
where $\w_i$ is the local disturbance. These disturbances are modeled by $S_1=[0~1;\,0~0]$, $F_1=[1~0]$, $S_2=0$,\,$F_2=1$, $S_3=[0~1;\,-1~0]$, $F_3=[1~0]$, $S_4=[0~2;\,-2~0]$, $F_4=[1~0]$. The local cost functions are chosen as $f_1(y_1)=\frac{1}{2}(y_1+2)^2$, $f_2(y_2) = y_2^2 \ln(1 + y_2^2) + (y_2+1)^2$, $f_3(y_3)=\ln(e^{-0.1y_3}+e^{0.3y_3})+y_3^2$, 
$f_4(y_4)=\frac{y_4^2}{25\sqrt{y_4^2 + 1}}+(y_4-3)^2.$ Set $d_i=i$ and graph $\mathcal{G}$ as in Fig.~\ref{fig:graph}. The optimal solution is $y^*\approx \mbox{col}(2.2,\,0.7,\,2.0,\,5.1)$.

Letting $x_{i1}=y_i,\, x_{i2}=\dot{y}_i$ gives
\begin{align*}
\dot{x}_{i1}=x_{i2},\quad \dot{x}_{i2}=-f_ix_{i1}-g_ix_{i2}+u_i+F_id_i, \quad y_i=x_{i1}.
\end{align*}
It can be verified that Assumptions \ref{ass:graph}--\ref{ass:regulator} are all satisfied. Thus, the associated distributed optimal steady-state regulation problem can be solved by the controls given in last sections.

\begin{figure}
	\centering
	\includegraphics[width=0.84\textwidth]{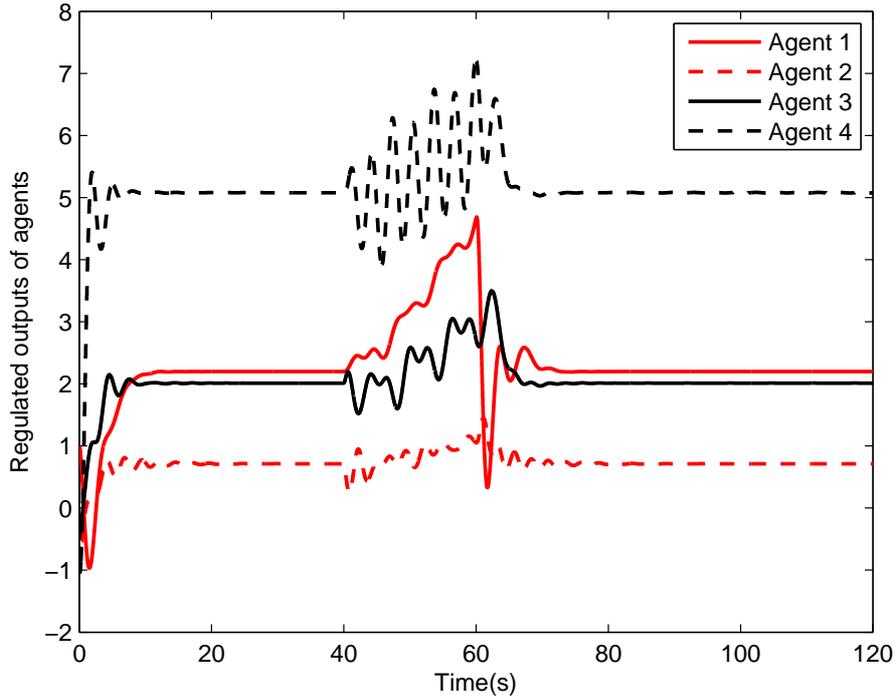}
	\caption{Profiles of all outputs under state-feedback control \eqref{ctrl:main-online}.}\label{fig:simu-3}
\end{figure}

For simulations, we assume only real-time gradients are available. The system parameters are taken as $f_1=1,\, g_1=1,\,f_2=0,\, g_2=1$,\, $f_3=1,\,g_3=0$,\, $f_4=0,\,g_4=0$. Then, we can solve the regulator equations in Assumption \ref{ass:regulator} and choose proper gain matrices for the state feedback control \eqref{ctrl:main-online}. To verify disturbance rejection performance, we first simulate the disturbance-free dynamics during $t= 0\sim 40\,{\rm s}$, and then the disturbed case after $t=40\,{\rm s}$. After the disturbance rejection part works at $t=60\,{\rm s}$, the optimal steady-state is quickly recovered. The simulation result is showed in Fig.~\ref{fig:simu-3}.  Satisfactory regulating performance is observed.

\section{Conclusions}\label{sec:con}

A distributed optimal steady-state regulation problem was formulated for general linear agents as a combination of asymptotic steady-state regulation and distributed resource allocation. By an embedded control scheme, we introduced an optimal signal generator and divided this problem into two simpler subproblems. Both state and output feedback cases were investigated and solved with disturbance rejection. A nontrivial extension using only real-time gradient information was also presented by high-gain control techniques. Future works will include nonlinear agents with nonsmooth cost functions.

\bibliographystyle{IEEEtran}
\bibliography{opt-regulation}

\end{document}